\documentclass{birkmult}
\newcommand{\NN}{{\mathbb{N}}}

\newcommand{\eps}{\varepsilon}
\newcommand{\po}{\partial\Omega}

\newcommand{\bp}{\noindent {\it Proof}.\,\,}
\newcommand{\ep}{\hfill$\Box$ \vskip 0.08in}

\def\ring{\mathaccent"0017 }

 \newtheorem{proposition}{Proposition}[section]
\newtheorem{theorem}[proposition]{Theorem}

\numberwithin{equation}{section}

\begin{document}
%
%
%
%
%
%
%
%
%
\title[Polyharmonic Green function]{Pointwise estimates for the polyharmonic Green function in general domains}

\author[Svitlana Mayboroda]{Svitlana Mayboroda}

\address{%
Department of Mathematics\\
Purdue University\\
150 N. University Street\\
West Lafayette, IN 47907-2067}

\email{svitlana@math.purdue.edu}

\thanks{ The first author is partially supported by the NSF grant DMS 0758500.}

\author[Vladimir Maz'ya]{Vladimir Maz'ya}
\address{%
Department of Mathematical Sciences, M\&O\\
  Building, University of Liverpool, Liverpool L69 3BX, UK\\
  Department of Mathematics, Link\"oping University, \\
  SE-581 83 Link\"oping, Sweden}
  
 \email{vlmaz@mai.liu.se}




\maketitle

\section{Introduction}
\setcounter{equation}{0}

Let $\Omega$ be a arbitrary bounded domain in ${\mathbb{R}}^n$.  The Green
function for the polyharmonic equation is a function $G:\Omega\times\Omega\to{\mathbb{R}}$ which for every
fixed $y\in\Omega$ solves the equation
\begin{equation}\label{eq4.1}
(-\Delta_x)^m G(x,y)=\delta(x-y), \qquad x\in\Omega,
\end{equation}

\noindent in the space $\ring W^{m,2}(\Omega)$, a
completion of $C_0^\infty(\Omega)$ in the norm given by $\|u\|_{\ring
W^{m,2}(\Omega)}=\|\nabla^m u\|_{L^2(\Omega)}$. The case $m=2$ corresponds to the biharmonic equation, and respectively, \eqref{eq4.1} gives rise to the biharmonic Green function.

In dimension two the biharmonic Green function can be interpreted as a deflection of a thin clamped plate under a point load. Numerous applications in structural engineering, emerging from this fact, have stimulated considerable interest to the biharmonic equation and its Green function as early as in the beginning of 20th century. In 1908 Hadamard has published a volume devoted to properties of the solutions to the biharmonic equation \cite{Hadamard}, where, in particular, he conjectured that the corresponding Green function must be positive, at least, in convex domains. However, several  counterexamples to Hadamard's conjecture have been found later on (\cite{Duffin1}, \cite{Duffin2}, \cite{Garabedian}, \cite{Loewner}, \cite{Szego}, \cite{CoffmanGrover}, \cite{KKM}) and it was proved that the biharmonic Green function may change sign even in a smooth convex domain, in a sufficiently eccentric ellipse (\cite{Garabedian}, \cite{CoffmanGrover}). Moreover, in a rectangle the first eigenfunction of the biharmonic operator has infinitely many changes of sign near each of the vertices (\cite{Coffman}, \cite{KKM}).

During the past century, the biharmonic and more generally, the polyharmonic Green function has been thoroughly studied, and a variety of upper estimates has been obtained. In particular, we would like to point out the results in smooth domains \cite{DASw}, \cite{Kras}, \cite{Sol1}, \cite{Sol2}, in conical domains \cite{MazPlam}, \cite{KMRElliptic}, and in polyhedra  \cite{MazRoss}. 

{\it The objective of the present paper is to establish sharp estimates on the polyharmonic Green function and its derivatives without any geometric assumptions, in an arbitrary bounded open set. }

For example, we show that, whenever the dimension $n\in [3,2m+1]\cap\NN$ is odd, the regular part of the Green function admits the estimate
\begin{equation}\label{int1}
|\nabla^{m-\frac n2+\frac 12}_x\nabla^{m-\frac n2+\frac 12}_y(G(x,y)-\Gamma(x-y))|\leq {{\frac{C}{\max\{d(x),d(y),|x-y|\}}}}, \quad x,y\in\Omega,
\end{equation}

\noindent where $\Gamma(x)=C_{m,n}|x|^{2m-n}$, $x\in\Omega$, is a fundamental solution for the polyharmonic operator,  $d(x)$ is the distance from $x\in\Omega$ to
$\po$ and the constant $C$ depends on $n$ and $m$ only. Hence, in particular, 
\begin{equation}\label{int2}
|\nabla^{m-\frac n2+\frac 12}_x\nabla^{m-\frac n2+\frac 12}_yG(x,y)|\leq {{\frac{C}{|x-y|}}},\qquad x,y\in\Omega,
\end{equation}

\noindent and similar results are established for the lower order derivatives. 

Furthermore, the estimates on the Green function allow us to derive optimal bounds for the solution $u$ of the Dirichlet boundary value problem 
\begin{equation}\label{int3}
(-\Delta)^m u=\sum_{|\alpha|\leq m-\frac n2+\frac 12}c_\alpha \partial^{\alpha}f_{\alpha}, \quad u\in \ring W^{m,2}(\Omega).
\end{equation}

\noindent  Specifically,
\begin{eqnarray}\label{int4}|\nabla^{m-\frac n2+\frac 12} u(x)|&\leq & C  \sum_{|\alpha|\leq  m-\frac n2+\frac 12}
\int_\Omega d(y)^{m-\frac n2+\frac 12-|\alpha|}\,\frac{|f_\alpha(y)|}{|x-y|}\,dy,\,\,  x\in\Omega,
\end{eqnarray}

\noindent whenever the integrals on the right-hand side of \eqref{int4} are finite. 
In particular, there exists a constant $C_{\Omega}>0$ depending on $m$, $n$ and the domain $\Omega$ such that 
\begin{eqnarray}\label{int5}\|\nabla^{m-\frac n2+\frac 12} u\|_{L^\infty(\Omega)}&\leq & C_\Omega    \sum_{|\alpha|\leq  m-\frac n2+\frac 12}
\| d(\cdot)^{m-\frac n2-\frac 12-|\alpha|}f_{\alpha}\|_{L^p(\Omega)},\,\, 
\end{eqnarray}

\noindent for $p>\frac{n}{n-1}$.

The bounds above are sharp,  in the sense that the solution of the polyharmonic equation in an arbitrary domain generally does not exhibit more regularity.
Indeed, assume that $n\in [3,2m+1]\cap\NN$ is odd and let $\Omega\subset{\mathbb{R}}^n$ be the punctured unit ball
$B_1\setminus\{O\}$, where $B_r=\{x\in{\mathbb{R}}^n:\,|x|<r\}$. Consider
a function $\eta\in C_0^\infty(B_{1/2})$ such that $\eta=1$ on
$B_{1/4}$. Then let
\begin{equation}\label{int6}
u(x):=\eta(x)\,\partial_x^{m-\frac n2-\frac 12}\Gamma(x)=C\eta(x)\,\partial_x^{m-\frac n2-\frac 12}(|x|^{2m-n}),\qquad x\in B_1\setminus\{O\},
\end{equation}

\noindent where $\partial_x$ stands for a derivative in the direction of $x_i$ for some $i=1,...,n$.
It is straightforward to check that $u\in \ring W^{m,2}(\Omega)$ and $(-\Delta)^m u\in
C_0^\infty (\Omega)$. While $\nabla^{m-\frac n2+\frac 12} u$ is bounded, the derivatives of the order $m-\frac n2+\frac 32$ are not, and moreover, $\nabla^{m-\frac n2+\frac 12} u$
is not
continuous at the origin. Therefore, the estimates \eqref{int4} are optimal in general domains.

We also derive full analogues of \eqref{int1}, \eqref{int2}, \eqref{int4}, \eqref{int5} and accompanying lower order bounds in even dimensions. In that case, the optimal regularity turns out to be of the order $m-\frac n2$.

Finally, we would like to mention that 
the Green function estimates in this paper generalize the earlier developments in \cite{MayMazBiharm}, where the biharmonic Green function was treated, and \cite{MazyaCh}, where the pointwise estimates on polyharmonic Green function have been established in dimensions $2m+1$ and $2m+2$ for  $m>2$ and dimensions $5,6,7$ for $m=2$.

\section{Preliminaries}
\setcounter{equation}{0}

The Green function estimates in the present paper are based, in particular, on the recent results
for locally polyharmonic functions that will appear in \cite{MayMazPolyharm}. We record them below without the proof. 
Here and throughout the paper $B_r(Q)$ and $S_r(Q)$ denote,
respectively, the ball and the sphere with radius $r$ centered at
$Q$ and $C_{r,R}(Q)=B_R(Q)\setminus\overline{B_r(Q)}$. When the
center is at the origin, we write $B_r$ in place of $B_r(O)$, and
similarly $S_r:=S_r(O)$ and $C_{r,R}:=C_{r,R}(O)$. Also, $\nabla^mu$
stands for a vector of all derivatives of $u$ of the order $m$.

\begin{proposition}\label{c3.3}
Let $\Omega$ be a bounded domain in ${\mathbb{R}}^n$, $2\leq n \leq 2m+1$,
$Q\in{\mathbb{R}}^n\setminus\Omega$, and $R>0$. Suppose
\begin{equation}\label{eq3.3}
(-\Delta)^m u=f \,\,{\mbox{in}}\,\,\Omega, \quad f\in
C_0^{\infty}(\Omega\setminus B_{4R}(Q)),\quad u\in \ring
W^{m,2}(\Omega).
\end{equation}

\noindent Then
\begin{equation}\label{eq3.4}
\frac{1}{\rho^{2\lambda +n-1}}\int_{S_{\rho}(Q)\cap\Omega}|u(x)|^2\,d\sigma_x
\leq \frac{C}{R^{2\lambda +n}} \int_{C_{R,4R}(Q)\cap\Omega} |u(x)|^2\,dx\quad
{\mbox{ for every}}\quad \rho<R,
\end{equation}

\noindent where $C$ is a  constant depending on $m$ and $n$ only, and 
\begin{equation}\label{eq3.4-1}
\lambda=m-n/2+1/2 \,\,\mbox{when $n$ is odd},\qquad \quad \lambda=m-n/2 \,\,\mbox{when $n$ is even}.
\end{equation}

Moreover, for every $x\in B_{R/4}(Q)\cap\Omega$
\begin{equation}\label{eq3.14}
|\nabla^{i} u(x)|^2\leq C \frac{|x-Q|^{2\lambda-2i}}{R^{n+2\lambda}}\int_{C_{R/4,4R}(Q)\cap\Omega}
|u(y)|^2\,dy,\qquad 0\leq i\leq \lambda
\end{equation}

\noindent where $\lambda$ is given by \eqref{eq3.4-1}.
\end{proposition}

In addition, using the Kelvin transform, estimates near the origin for solutions of elliptic equations can be translated into estimates at infinity. In particular, Proposition~\ref{c3.3} leads to the following result (also proved in \cite{MayMazPolyharm}).

\begin{proposition}\label{p3.4} Let $\Omega$ be a bounded domain in ${\mathbb{R}}^n$, $2\leq n \leq 2m+1$,
$Q\in{\mathbb{R}}^n\setminus\Omega$, $r>0$ and assume that
\begin{equation}\label{eq3.22}
(-\Delta)^m u=f \,\,{\mbox{in}}\,\,\Omega, \quad f\in
C_0^{\infty}(B_{r/4}(Q)\cap\Omega),\quad u\in \ring W^{m,2}(\Omega).
\end{equation}

\noindent Then
\begin{equation}\label{eq3.23}
\rho^{2\lambda+n+1-4m}\int_{S_{\rho}(Q)\cap\Omega}|u(x)|^2\,d\sigma_x
\leq C\, r^{2\lambda+n-4m}\int_{ C_{r/4,r}(Q)\cap\Omega} |u(x)|^2\,dx,
\end{equation}

\noindent for any $\rho>r$ and $\lambda$ given by \eqref{eq3.4-1}.

Furthermore, for any $x\in\Omega\setminus B_{4r}(Q)$
\begin{equation}\label{eq3.24}
|\nabla^i u(x)|^2\leq
C\,\frac{r^{2\lambda+n-4m}}{|x-Q|^{2\lambda+2n-4m+2i}}\int_{C_{r/4,4r}(Q)\cap\Omega} |u(y)|^2\,dy,\qquad 0\leq i\leq \lambda.
\end{equation}
\end{proposition}

\section{Estimates for the Green function}\label{SGrMain}
\setcounter{equation}{0}

Following \cite{PolyharmBook} we point out that the fundamental solution for the $m$-Laplacian is a linear combination of the characteristic singular solution (defined below) and any $m$-harmonic function in ${\mathbb{R}}^n$. The characteristic singular solution is
\begin{eqnarray}\label{chfundsol1}
&C_{m,n}|x|^{2m-n},\quad &\mbox{if $n$ is odd, or if $n$ is even with $n\geq 2m+2$},\\[4pt]
\label{chfundsol2}
&C_{m,n}|x|^{2m-n} \log |x|,\quad &\mbox{if $n$ is even with $n\leq 2m$}.
\end{eqnarray}

\noindent The exact expressions for constants  $C_{m,n}$ can be found in \cite{PolyharmBook}, p.8.  For the purposes of this paper we will use the fundamental solution given by 
\begin{equation}\label{fundsol1}
\Gamma(x)=C_{m,n}  \left\{\begin{array}{l}
|x|^{2m-n},\quad \mbox{if $n$ is odd},\\[4pt]
|x|^{2m-n}\log \frac{{\rm diam}\,\Omega}{|x|},\quad \mbox{if $n$ is even and $n\leq 2m$},\\[4pt]
|x|^{2m-n},\quad \mbox{if $n$ is even and $n\geq 2m+2$}.
\end{array}
\right.
\end{equation}

\begin{theorem}\label{p4.1} Let $\Omega\subset {\mathbb{R}}^n$ be a bounded
domain. Then there exist  constants $C$, $C'$ depending on $m$ and $n$ only such that for
every $x,y\in\Omega$ the following estimates hold. If $n\in [3,2m+1] \cap\NN$ is odd then
\begin{eqnarray}\label{Gr1}\nonumber
|\nabla_x^i\nabla_y^j G(x,y)| &\leq & C \min\left\{1,\Bigg(\frac{d(x)}{|x-y|}\Bigg)^{m-\frac n2+\frac 12-i}, \Bigg(\frac{d(y)}{|x-y|}\Bigg)^{m-\frac n2+\frac 12-j}\right\}\\[4pt]
&&\quad \times \frac{1}{|x-y|^{n-2m+i+j}},
\end{eqnarray}

\noindent whenever $0\leq i,j\leq m-\frac n2+\frac 12$ are such that $i+j\geq 2m-n$, and 
\begin{eqnarray}\label{Gr2}\nonumber
|\nabla_x^i\nabla_y^j G(x,y)| &\leq & C \min\left\{1,\Bigg(\frac{d(x)}{|x-y|}\Bigg)^{m-\frac n2+\frac 12-i}, \Bigg(\frac{d(y)}{|x-y|}\Bigg)^{m-\frac n2+\frac 12-j}\right\}\times\\[4pt] && \,\times\frac{1}{|x-y|^{n-2m+i+j}}\min\left\{\frac{|x-y|}{d(x)},\frac{|x-y|}{d(y)},1\right\}^{n-2m+i+j},
\end{eqnarray}

\noindent if $0\leq i,j\leq m-\frac n2+\frac 12$ are such that $i+j\leq 2m-n$. 

If $n\in [2,2m]\cap\NN$ is even, then 
\begin{eqnarray}\label{Gr3}\nonumber
|\nabla_x^i\nabla_y^j G(x,y)| &\leq & C \min\left\{1,\Bigg(\frac{d(x)}{|x-y|}\Bigg)^{m-\frac n2-i}, \Bigg(\frac{d(y)}{|x-y|}\Bigg)^{m-\frac n2-j}\right\}\times\\[4pt] \nonumber && \,\times\frac{1}{|x-y|^{n-2m+i+j}}\min\left\{\frac{|x-y|}{d(x)},\frac{|x-y|}{d(y)},1\right\}^{n-2m+i+j} \times\\[4pt] && \,\times\log \left(1+\frac{\min\{d(x),d(y)\}}{|x-y|}\right),
\end{eqnarray}

\noindent for all $0\leq i,j\leq m-\frac n2$. 

Furthermore, the estimates on the regular part of the Green function $S(x,y)=G(x,y)-\Gamma(x-y)$, $x,y\in\Omega$, are as follows.  If $n\in [3,2m+1] \cap\NN$ is odd then
\begin{eqnarray}\label{Gr1s}
&&|\nabla_x^i\nabla_y^j S(x,y)|\leq \frac{C}{\max\{d(x),d(y),|x-y|\}^{n-2m+i+j}},
\end{eqnarray}
\noindent whenever $0\leq i,j\leq m-\frac n2+\frac 12$ are such that $i+j\geq 2m-n$, and 
\begin{eqnarray}\label{Gr2s}
|\nabla_x^i\nabla_y^j S(x,y)| \leq  \frac{C}{|x-y|^{n-2m+i+j}}\,\min\left\{\frac{|x-y|}{d(x)},\frac{|x-y|}{d(y)},1\right\}^{n-2m+i+j},
\end{eqnarray}

\noindent if $0\leq i,j\leq m-\frac n2+\frac 12$ are such that $i+j\leq 2m-n$. 

If $n\in [2,2m]\cap\NN$ is even, then 
\begin{eqnarray}\label{Gr3s}\nonumber
|\nabla_x^i\nabla_y^j S(x,y)| &\leq &\frac{C}{|x-y|^{n-2m+i+j}}\min\left\{\frac{|x-y|}{d(x)},\frac{|x-y|}{d(y)},1\right\}^{n-2m+i+j}\times\\[4pt] && \,\times \log \left(1+\frac{{\rm diam}\,\Omega}{\max\{d(x),d(y),|x-y|\}}\right),
\end{eqnarray}

\noindent for all $0\leq i,j\leq m-\frac n2$. 
\end{theorem}


\bp Let us start with some auxiliary calculations. Let $\alpha$ be a multi-index of length less than or equal to $\lambda$, where 
$\lambda$ is given by \eqref{eq3.4-1}. Then $\partial^{\alpha}_y\,\Gamma(x-y)$ can be written as 
\begin{equation}\label{eqGr1}
\partial^{\alpha}_y\,\Gamma(x-y)=P^{\alpha}(x-y)\log \frac{{\rm diam}\,\Omega}{|x-y|}+Q^{\alpha}(x-y).
\end{equation}

\noindent When the dimension is odd, $P^{\alpha}\equiv 0$. If the dimension is even (and less than or equal to $2m$ by the assumptions of the theorem) then $P^{\alpha}$ is a homogeneous polynomial of order $2m-n-|\alpha|$ as long as $|\alpha|\leq 2m-n$. In any case, $Q^{\alpha}$ is a homogeneous function of order $2m-n-|\alpha|$.

Consider a
function $\eta$ such that
\begin{equation}\label{eq4.4}
\eta \in C_0^\infty(B_{1/2})\quad\mbox{and}\quad \eta=1
\quad\mbox{in}\quad B_{1/4},
\end{equation}

\noindent and define 
\begin{equation}\label{eq4.5}
{\mathcal R}_{\alpha} (x,y):=\partial^{\alpha}_y\,
G(x,y)-\eta\left(\frac{x-y}{d(y)}\right)\left(P^{\alpha}(x-y)\log \frac{d(y)}{|x-y|}+Q^{\alpha}(x-y)\right),
\end{equation}

\noindent  for $x,y\in\Omega.$
Also, let us denote
\begin{eqnarray}\label{eq4.6}\nonumber
&&\hskip -1cm f_{\alpha}(x,y):= (-\Delta_x)^m {\mathcal R}_{\alpha}
(x,y)\\[4pt]&&\hskip -1cm\qquad =-\left[(-\Delta_x)^m,\eta\left(\frac{x-y}{d(y)}\right)
\right]\left(P^{\alpha}(x-y)\log \frac{d(y)}{|x-y|}+Q^{\alpha}(x-y)\right).
\end{eqnarray}

\noindent It is not hard to see that for every $\alpha$ as above
\begin{equation}\label{eq4.7}
f_{\alpha}(\cdot,y)\in C_0^\infty(C_{d(y)/4,
d(y)/2}(y))\qquad\mbox{and}\qquad |f_\alpha(x,y)|\leq Cd(y)^{-n-|\alpha|}, \quad
x,y\in\Omega.
\end{equation}

\noindent Then for every fixed $y\in\Omega$ the function $x\mapsto
{\mathcal R}_{\alpha}(x,y)$ is a solution of the boundary value problem
\begin{equation}\label{eq4.8}
(-\Delta_x)^m {\mathcal R}_\alpha (x,y)= f_\alpha(x,y) \,\,{\mbox{in}}\,\,\Omega,
\quad f_\alpha(\cdot,y)\in C_0^{\infty}(\Omega),\quad {\mathcal
R}_\alpha(\cdot,y)\in \ring W^{m,2}(\Omega),
\end{equation}

\noindent   so that
\begin{equation}\label{eq4.9}
\left\|\nabla_x^m{\mathcal
R}_\alpha(\cdot,y)\right\|_{L^2(\Omega)}=\left\|{\mathcal
R}_\alpha(\cdot,y)\right\|_{W^{m,2}(\Omega)}\leq
C\|f_\alpha(\cdot,y)\|_{W^{-m,2}(\Omega)},\qquad 0\leq|\alpha|\leq \lambda.
\end{equation}

\noindent Here $W^{-m,2}(\Omega)$ stands for the Banach space dual
of $\ring W^{m,2}(\Omega)$, i.e.
\begin{equation}\label{eq4.10}
\|f_\alpha(\cdot,y)\|_{W^{-m,2}(\Omega)}=\sup_{v\in \ring
W^{m,2}(\Omega):\,\|v\|_{\ring
W^{m,2}(\Omega)}=1}\int_{\Omega}f_\alpha(x,y)v(x)\,dx.
\end{equation}

 Recall that by Hardy's inequality
\begin{equation}\label{eq4.11}
\left\|\frac{v}{|\cdot-\,Q|^{m}}\right\|_{L^2(\Omega)}\leq C
\left\|\nabla^m v\right\|_{L^2(\Omega)}\quad\mbox{for every}\quad
v\in\ring W^{m,2}(\Omega),\quad Q\in\po.
\end{equation}

\noindent Then for some $y_0\in\po$ such that $|y-y_0|=d(y)$ and any $v$ in \eqref{eq4.10}
\begin{eqnarray}\label{eq4.12}\nonumber
&&\int_{\Omega}f_\alpha(x,y)v(x)\,dx\leq C
\left\|\frac{v}{|\cdot-y_0|^{m}}\right\|_{L^2(\Omega)}\left\|f_\alpha(\cdot,y)
|\cdot-y_0|^{m}\right\|_{L^2(\Omega)} \\[4pt]
&&\qquad \leq C d(y)^{m} \left\|\nabla^m v\right\|_{L^2(\Omega)}
\|f_\alpha(\cdot,y) \|_{L^2(C_{d(y)/4, d(y)/2}(y))},
\end{eqnarray}

\noindent and therefore, by (\ref{eq4.7})
\begin{equation}\label{eq4.13}
\left\|\nabla_x^m{\mathcal R}_\alpha(\cdot,y)\right\|_{L^2(\Omega)}\leq C
d(y)^{m-|\alpha|-n/2}.
\end{equation}

Now we split the discussion into a few cases. 

\vskip 0.08in \noindent {\bf Case I: $|x-y|\geq N d(y)$ or $|x-y|\geq N d(x)$} for some large $N$ to be specified later.

Let us  first assume that
$|x-y|\geq N d(y)$. As
before, we denote by $y_0$ some point on the boundary such that
$|y-y_0|=d(y)$. Then by (\ref{eq4.7})--\eqref{eq4.8} the function $x\mapsto {\mathcal
R}_\alpha(x,y)$ is $m$-harmonic in $\Omega\setminus B_{3d(y)/2}(y_0)$.
Hence, by Proposition~\ref{p3.4} with $r=6d(y)$
\begin{equation}\label{eq4.14}
|\nabla_x^i {\mathcal R}_{\alpha}(x,y)|^2\leq
C\,\frac{d(y)^{2\lambda+n-4m}}{|x-y_0|^{2\lambda+2n-4m+2i}}\,\int_{C_{3d(y)/2,24d(y)}(y_0)}|{\mathcal
R}_\alpha(z,y)|^2\,dz,
\end{equation}

\noindent provided that $0\leq i\leq \lambda$ and $|x-y|\geq 4r+d(y)$, i.e $N\geq 25$. The
right-hand side of (\ref{eq4.14}) is bounded by
\begin{eqnarray}\label{eq4.15}
&& C\,\frac{d(y)^{2\lambda+n-2m}}{|x-y_0|^{2\lambda+2n-4m+2i}}\,\int_{C_{3d(y)/2,24d(y)}(y_0)}\frac{|{\mathcal
R}_\alpha(z,y)|^2}{|z-y_0|^{2m}}\,dz\nonumber\\[4pt]
&&\qquad \leq C\,\frac{d(y)^{2\lambda+n-2m}}{|x-y_0|^{2\lambda+2n-4m+2i}}\,\int_{\Omega}|\nabla_z^m{\mathcal R}_\alpha(z,y)|^2\,dz
\nonumber\\[4pt]
&&\qquad \leq C\,\frac{d(y)^{2\lambda-2|\alpha|}}{|x-y|^{2\lambda+2n-4m+2i}},
\end{eqnarray}

\noindent by Hardy's inequality and (\ref{eq4.13}). Therefore, 
\begin{equation}\label{eq4.15-1}
|\nabla_x^i {\mathcal R}_{\alpha}(x,y)|^2\leq
C\,\frac{d(y)^{2\lambda-2|\alpha|}}{|x-y|^{2\lambda+2n-4m+2i}}, \quad \mbox{when}\quad |x-y|\geq N d(y).\quad 0\leq i,|\alpha|\leq \lambda.
\end{equation}

\noindent Since for $N\geq 25$ the condition $|x-y|\geq N d(y)$ guarantees that $\eta\left(\frac{x-y}{d(y)}\right)=0$ and hence, ${\mathcal R}_{\alpha}(x,y)=\partial_y^\alpha G(x,y)$ when $|x-y|\geq Nd(y)$, the estimate \eqref{eq4.15-1} with $j:=|\alpha|$ implies
\begin{equation}\label{eq4.15-2}
|\nabla_x^i \nabla_y^{j}G(x,y)|^2\leq
C\,\frac{d(y)^{2\lambda-2j}}{|x-y|^{2\lambda+2n-4m+2i}}, \quad \mbox{when}\quad |x-y|\geq N d(y),\quad 0\leq i,j\leq \lambda.
\end{equation}

\noindent Also, by the symmetry of the Green function we automatically deduce that
\begin{equation}\label{eq4.15-2s}
|\nabla_x^i \nabla_y^{j}G(x,y)|^2\leq
C\,\frac{d(x)^{2\lambda-2i}}{|x-y|^{2\lambda+2n-4m+2j}}, \quad \mbox{when}\quad |x-y|\geq N d(x),\quad 0\leq i,j\leq \lambda.
\end{equation}

\noindent In particular, \eqref{eq4.15-2} and \eqref{eq4.15-2s} combined give the estimate
\begin{equation}\label{eq4.15-3}
|\nabla_x^i \nabla_y^{j}G(x,y)|\leq
\frac{C}{|x-y|^{n-2m+i+j}}, \quad \mbox{when}\quad |x-y|\geq N \min\{d(x),d(y)\},
\end{equation}

\noindent for $0\leq i,j\leq \lambda$.

 Now further consider several cases. If $n$ is odd, then 
\begin{equation}\label{eq4.16}
|\nabla_x^i\nabla_y^j \Gamma(x-y)|\leq \frac{C}{|x-y|^{n-2m+i+j}}\qquad\mbox{for
all}\qquad x,y\in\Omega, \quad i,j\geq 0,
\end{equation}

\noindent while if $n$ is even, then 
\begin{equation}\label{eq4.16-1}
|\nabla_x^i\nabla_y^j \Gamma(x-y)|\leq C_1\,|x-y|^{-n+2m-i-j}\log \frac{{\rm diam}\,{(\Omega)}}{|x-y|}+ C_2\, |x-y|^{-n+2m-i-j},
\end{equation}

\noindent for all $x,y\in\Omega$ and $0\leq i+j\leq 2m-n$. 

Combining this  with  (\ref{eq4.15-3}) we
deduce that for $n\leq 2m+1$ odd
\begin{equation}\label{eq4.16s}
|\nabla_x^i\nabla_y^j S(x,y)|\leq \frac{C}{|x-y|^{n-2m+i+j}}\quad \mbox{when}\quad |x-y|\geq N \min\{d(x),d(y)\},\end{equation}

\noindent while if $n$ is even, then 
\begin{equation}\label{eq4.16-1s}
|\nabla_x^i\nabla_y^j S(x,y)|\leq C\,|x-y|^{-n+2m-i-j}\left(C'+\log \frac{{\rm diam}\,{(\Omega)}}{|x-y|}\right)\quad 
\end{equation}

\noindent provided that $|x-y|\geq N \min\{d(x),d(y)\}$ and $0\leq i,j\leq \lambda$.


\vskip 0.08in \noindent {\bf Case II: $|x-y|\leq N^{-1} d(y)$ or $|x-y|\leq N^{-1} d(x)$}. 

Assume that $|x-y|\leq N^{-1}d(y)$. For such $x$ we have
$\eta\bigl(\frac{x-y}{d(y)}\bigr)= 1$ and therefore
\begin{equation}\label{eq4.19}
{\mathcal R}_{\alpha} (x,y)=\partial^{\alpha}_y\,
G(x,y)-P^{\alpha}(x-y)\log \frac{d(y)}{|x-y|}-Q^{\alpha}(x-y).
\end{equation}

\noindent Hence, if $n$ is odd, 
\begin{equation}\label{eq4.19-1}
{\mathcal R}_{\alpha} (x,y)=\partial^{\alpha}_y\left(
G(x,y)-
\Gamma(x-y)\right),\quad \mbox{when}\quad |x-y|\leq N^{-1}d(y),
\end{equation}
\noindent and if $n$ is even, 
\begin{equation}\label{eq4.19-2}
{\mathcal R}_{\alpha} (x,y)=\partial^{\alpha}_y\left(
G(x,y)-
\Gamma(x-y)\right)+P^\alpha(x-y)\log \frac{{\rm diam}\,\Omega}{d(y)},
\end{equation}

\noindent when $|x-y|\leq N^{-1}d(y)$.
By the interior estimates for solutions of elliptic
equations
\begin{equation}\label{eq4.20}
|\nabla_x^i{\mathcal R}_\alpha(x,y)|^2\leq
\frac{C}{d(y)^{n+2i}}\int_{B_{d(y)/8}(x)}|{\mathcal R}_\alpha(z,y)|^2\,dz,\quad\mbox{ for any $i\leq m$},
\end{equation}

\noindent since the function ${\mathcal R}_\alpha$ is $m$-harmonic in
$B_{d(y)/4}(y)\supset B_{d(y)/8}(x)$. Now we  bound the expression
above by
\begin{eqnarray}\label{eq4.21}
\nonumber&&\frac{C}{d(y)^{n+2i-2m}}\int_{B_{d(y)/4}(y)}\frac{|{\mathcal
R}_\alpha(z,y)|^2}{|z-y_0|^{2m}}\,dz \leq \frac{C}{d(y)^{n+2i-2m}}\left\|\nabla_x^m
{\mathcal R}(\cdot,y)\right\|_{L^2(\Omega)}^2\\[4pt]&&\qquad \leq \frac{C}{d(y)^{2n-4m+2i+2|\alpha|}},
\end{eqnarray}

\noindent with $0\leq |\alpha|\leq \lambda$. 

Let us now focus on the case of $n$ odd. It follows from \eqref{eq4.19-1} and \eqref{eq4.20} -- \eqref{eq4.21} that   
\begin{equation}\label{eq4.20s}
|\nabla_x^i\nabla_y^jS(x,y)|\leq
\frac{C}{d(y)^{n-2m+i+j}}, \quad 0\leq i\leq m, \quad 0\leq j\leq \lambda, \quad |x-y|\leq N^{-1}d(y),
\end{equation}

\noindent and hence, by symmetry,
\begin{equation}\label{eq4.20s-1}
|\nabla_x^i\nabla_y^jS(x,y)|\leq
\frac{C}{d(x)^{n-2m+i+j}}, \quad 0\leq i\leq\lambda, \quad 0\leq j\leq m, \quad |x-y|\leq N^{-1}d(x). 
\end{equation}

\noindent However,  we have
\begin{equation}\label{eq4.22}
|x-y|\leq N^{-1}d(y)\quad\Longrightarrow \quad (N-1)\, d(y)\leq N d(x)\leq (N+1)\, d(y),
\end{equation}

\noindent i.e. $d(y)\approx d(x)$ whenever $|x-y|$ is less than or equal to either $N^{-1}d(y)$ or $N^{-1}d(x)$. 
Therefore, when the dimension is odd, 
\begin{equation}\label{eq4.20s-2}
|\nabla_x^i\nabla_y^jS(x,y)|\leq
\frac{C}{\max\{d(x),d(y)\}^{n-2m+i+j}},  
\end{equation}

\noindent provided that $|x-y|\leq N^{-1}\max\{d(x),d(y)\},$ $0\leq i,j\leq\lambda$, $i+j\geq 2m-n$, and  
\begin{equation}\label{eq4.20s-3}
|\nabla_x^i\nabla_y^jS(x,y)|\leq
{C}\,{\min\{d(x),d(y)\}^{2m-n-i-j}},  
\end{equation}
\noindent for $|x-y|\leq N^{-1}\max\{d(x),d(y)\},$ $0\leq i,j\leq\lambda$, $i+j\leq 2m-n$.


As for the Green function itself, we then have for $ |x-y|\leq N^{-1}\max\{d(x),d(y)\}$
\begin{equation}\label{eq4.20g}
|\nabla_x^i\nabla_y^jG(x,y)|\leq
\frac{C}{|x-y|^{n-2m+i+j}}, \quad \mbox{if}\,\, i+j\geq 2m-n, \end{equation}

\noindent and 
\begin{equation}\label{eq4.20g-1}
|\nabla_x^i\nabla_y^jG(x,y)|\leq 
{C}\,{\min\{d(x),d(y)\}^{2m-n-i-j}}, \quad \mbox{if}\,\, i+j\leq 2m-n, 
\end{equation}

\noindent with $i,j$ such that $0\leq i,j\leq\lambda$. 

Similar considerations apply to the case when the dimension is even, leading to the following results: 
\begin{eqnarray}\label{eq4.20s-e}
&& |\nabla_x^i\nabla_y^jS(x,y)|\leq  
\frac{C}{d(y)^{n-2m+i+j}}\,\left( C'+\log \frac{{\rm diam}\,\Omega}{d(y)}\right), \,
\end{eqnarray}

\noindent for $ 0\leq i\leq m, \, 0\leq j\leq \lambda, \, |x-y|\leq N^{-1}d(y), $ and
\begin{eqnarray} \label{eq4.20s-1-e}
&& |\nabla_x^i\nabla_y^jS(x,y)|\leq 
\frac{C}{d(x)^{n-2m+i+j}}\,\left( C'+\log \frac{{\rm diam}\,\Omega}{d(x)}\right), \, 
\end{eqnarray}

\noindent for $0\leq i\leq\lambda, \, 0\leq j\leq m, \, |x-y|\leq N^{-1}d(x).$
In particular, in view of \eqref{eq4.22}, and the fact that $2m-n-i-j\geq 0$ whenever $0\leq i,j \leq \lambda$ and $n$ is even, we have
\begin{equation}\label{eq4.20s-2-e}
|\nabla_x^i\nabla_y^jS(x,y)|\leq
{C}\,{\min\{d(x),d(y)\}^{2m-n-i-j}}\left( C'+\log \frac{{\rm diam}\,\Omega}{{\max\{d(x),d(y)\}}}\right),  
\end{equation}

\noindent for $|x-y|\leq N^{-1}\max\{d(x),d(y)\}$, $0\leq i,j\leq\lambda$.


Passing to the Green function estimates, \eqref{eq4.19} and \eqref{eq4.20}--\eqref{eq4.21} lead to the bound 
\begin{eqnarray}\label{eq4.20g-e}
&& |\nabla_x^i\nabla_y^jG(x,y)|\leq  
{C}\,{d(y)^{2m-n-i-j}}\,\left( C'+\log \frac{d(y)}{|x-y|}\right), \, 
\end{eqnarray}

\noindent for $0\leq i,j\leq \lambda, \, |x-y|\leq N^{-1}d(y).$
Hence, by symmetry, 
\begin{eqnarray}\label{eq4.20g-1-e}
&& |\nabla_x^i\nabla_y^jG(x,y)|\leq  
{C}\,{d(x)^{2m-n-i-j}}\,\left( C'+\log \frac{d(x)}{|x-y|}\right), \, 
\end{eqnarray}

\noindent for $0\leq i,j\leq \lambda, \, |x-y|\leq N^{-1}d(x),$ and therefore, 
\begin{equation}\label{eq4.22-1}
|\nabla_x^i\nabla_y^jG(x,y)|\leq  
{C}\,{\min\{d(x),d(y)\}^{2m-n-i-j}}\,\left( C'+\log \frac{\min\{d(x),d(y)\}}{|x-y|}\right), \, 
\end{equation}

\noindent for all $0\leq i,j\leq \lambda,$ and  $ |x-y|\leq N^{-1}\max\{d(x),d(y)\}.$

Finally, it remains to consider  

\vskip 0.08in \noindent {\bf Case III: $|x-y|\approx d(y)\approx d(x)$},
 or more precisely, the situation when
\begin{equation}\label{eq4.25}
N^{-1}\,d(x)\leq |x-y|\leq Nd(x)\quad \mbox{and}\quad
N^{-1}\,d(y)\leq |x-y|\leq Nd(y).
\end{equation}

In this case we use the $m$-harmonicity of $x\mapsto G(x,y)$ in
$B_{d(x)/(2N)}(x)$.  Let $x_0\in\po$ be 
such that $|x-x_0|=d(x)$.
By the interior estimates, 
\begin{eqnarray}\label{eq4.26}\nonumber
&&|\nabla_x^i\nabla_y^{|\alpha|} G(x,y)|^2 \leq 
\frac{C}{d(x)^{n+2i}}\,\int_{B_{d(x)/(2N)}(x)}|\nabla_y^{|\alpha|}
G(z,y)|^2\,dz\\[4pt]
&&\quad \leq  \frac{C}{d(x)^{n+2i}}\,\int_{B_{d(x)/(2N)}(x)}|\nabla_y^{|\alpha|}
\Gamma(z-y)|^2\,dz\nonumber
\\[4pt]
&&\qquad \qquad
+
\frac{C}{d(x)^{n+2i-2m}}\,\int_{B_{2d(x)}(x_0)}\frac{|{\mathcal
R}_{\alpha}(z,y)|^2}{|z-x_0|^{2m}}\,dz\nonumber \\[4pt]
&&\quad \leq  \frac{C}{d(x)^{n+2i}}\,\int_{B_{d(x)/(2N)}(x)}|\nabla_y^{|\alpha|}
\Gamma(z-y)|^2\,dz\nonumber\\[4pt]
&&\qquad \qquad
+ \frac{C}{d(x)^{n+2i-2m}}\,\int_{\Omega}|\nabla_z^m {\mathcal
R}_{\alpha}(z,y)|^2\,dz\nonumber \\[4pt]
&&\quad \leq  \frac{C}{d(x)^{2n-4m+2i+2|\alpha|}}+\frac{C}{d(x)^{n-2m+2i}d(y)^{n-2m+2|\alpha|}},
\end{eqnarray}

\noindent provided that $0\leq i\leq m$, $0\leq |\alpha|\leq \lambda$ and $n$ is odd. The right-hand side of \eqref{eq4.26} also provided the estimate on derivatives of the Green function holds when $n$ is even, upon observing that

\begin{eqnarray}\label{eq4.26-e}\nonumber
&&
\frac{C}{d(x)^{n+2i}}\,\int_{B_{d(x)/(2N)}(x)}|\nabla_y^{|\alpha|}
G(z,y)|^2\,dz \leq \frac{C}{d(x)^{n+2i}}\,\int_{B_{d(x)/(2N)}(x)}|P^{\alpha}(z-y)|^2\,dz \\[4pt]
&&\qquad +\frac{C}{d(x)^{n+2i}}\,\int_{B_{d(x)/(2N)}(x)}|Q^{\alpha}(z-y)|^2\,dz\nonumber\\[4pt]
&&\qquad +
\frac{C}{d(x)^{n+2i-2m}}\,\int_{B_{2d(x)}(x_0)}\frac{|{\mathcal
R}_{\alpha}(z,y)|^2}{|z-x_0|^{2m}}\,dz,
\end{eqnarray}

\noindent since the absolute value of $\log \frac{|z-y|}{d(y)}$ is bounded by a constant for $z,x,y$ as in \eqref{eq4.26-e}, \eqref{eq4.25}. 

Hence, for $x,y$ satisfying \eqref{eq4.25} we have 
\begin{eqnarray}\label{eq4.26-1}\nonumber
|\nabla_x^i\nabla_y^jG(x,y)|&\leq &  \frac{C}{\min \{d(x),d(y),|x-y|\}^{n-2m+i+j}}
\\[4pt]
&\approx &\frac{C}{\max \{d(x),d(y),|x-y|\}^{n-2m+i+j}},
\end{eqnarray}

\noindent for $0\leq i,j\leq \lambda$.

When $n$ is odd, the same argument implies the following estimate on a regular part of Green function
\begin{eqnarray}\label{eq4.26-1s}\nonumber
|\nabla_x^i\nabla_y^jS(x,y)|&\leq & \frac{C}{\min \{d(x),d(y),|x-y|\}^{n-2m+i+j}}\\[4pt] & \approx & \frac{C}{\max \{d(x),d(y),|x-y|\}^{n-2m+i+j}},
\end{eqnarray}

\noindent for $0\leq i,j\leq \lambda$, and $x,y$ satisfying \eqref{eq4.25}. If $n$ is even, however, we are led to a bound
\begin{eqnarray}\label{eq4.26-2s}\nonumber
|\nabla_x^i\nabla_y^jS(x,y)|&\leq & \frac{C}{\min \{d(x),d(y),|x-y|\}^{n-2m+i+j}}\times
\\[4pt]\nonumber&&\quad\times
\left(C'+\log \frac{{\rm diam}\,\Omega}{\max \{d(x),d(y),|x-y|\}^{n-2m+i+j}}\right)\\[4pt]
&\approx & \frac{C}{\max \{d(x),d(y),|x-y|\}^{n-2m+i+j}}\nonumber\times\\[4pt]&&\quad\times\left(C'+\log \frac{{\rm diam}\,\Omega}{\max \{d(x),d(y),|x-y|\}^{n-2m+i+j}}\right)
\end{eqnarray}
\noindent for $0\leq i,j\leq \lambda$.

\vskip 0.08in

The final bounds for the Green function are a combination of estimates \eqref{eq4.15-2}, \eqref{eq4.15-2s}, \eqref{eq4.20g}, \eqref{eq4.20g-1}, \eqref{eq4.22-1}, \eqref{eq4.26-1}. It helps to observe that the regions of $(x,y)\in\Omega\times\Omega$ in \eqref{eq4.15-2},  \eqref{eq4.15-2s} are disjoint from those in \eqref{eq4.20g}, \eqref{eq4.20g-1}, \eqref{eq4.22-1}. The condition $|x-y|\leq N^{-1}\max\{d(x),d(y)\}$ excludes the possibility of $|x-y|\leq N^{-1}\max\{d(x),d(y)\}$. This is, in particular, due to \eqref{eq4.22}. Also, the bound \eqref{eq4.26-1} is the same as \eqref{eq4.15-2}, \eqref{eq4.15-2s}, \eqref{eq4.20g}, \eqref{eq4.20g-1} for the case when $d(x)$, $d(y)$ and $|x-y|$ are all comparable. Hence, it can be suitably absorbed. Finally, it is straightforward to check that 
\begin{equation}\label{eq4.26-1.1s}
C'+\log \frac{\min\{d(x),d(y)\}}{|x-y|}\approx  \log \left(1+\frac{\min\{d(x),d(y)\}}{|x-y|}\right)
\end{equation}
\noindent  for   $ |x-y|\leq N^{-1}\max\{d(x),d(y)\}.$

Analogously, the desired estimates on the regular part of the Green function can be drawn from 
\eqref{eq4.16s}, \eqref{eq4.16-1s}, \eqref{eq4.20s-2}, \eqref{eq4.20s-3}. \eqref{eq4.20s-2-e}, \eqref{eq4.26-1s}, \eqref{eq4.26-2s}.
\ep


\section{Applications: estimates on solutions of the Dirichlet problem} \setcounter{equation}{0}

Green function estimates proved in Section~\ref{SGrMain}
allow us to investigate the solutions of the Dirichlet problem
for the polyharmonic equation for a wide class of data. 

\begin{proposition}\label{pGrAp1} Let $\Omega\subset {\mathbb{R}}^n$ be a bounded
domain and assume that   $n\in [3,2m+1] \cap\NN$ is odd.  
Consider the
boundary value problem
\begin{equation}\label{eq4.27}
(-\Delta)^m u=\sum_{|\alpha|\leq m-\frac n2+\frac 12}c_\alpha \partial^{\alpha}f_{\alpha}\in W^{-m,2}(\Omega), \quad u\in \ring W^{m,2}(\Omega).
\end{equation}

\noindent Then the solution satisfies the estimate
\begin{eqnarray}\label{eq4.28}|\nabla^{m-\frac n2+\frac 12} u(x)|&\leq & C  \sum_{|\alpha|\leq  m-\frac n2+\frac 12}
\int_\Omega d(y)^{m-\frac n2+\frac 12-|\alpha|}\,\frac{|f_\alpha(y)|}{|x-y|}\,dy,\,\,  x\in\Omega,
\end{eqnarray}

\noindent whenever the integrals on the right-hand side of \eqref{eq4.28} are finite. The constant $C$ in \eqref{eq4.28} depends on $m$ and $n$ only.

In particular, there exists a constant $C_{\Omega}>0$ depending on $m$, $n$ and the domain $\Omega$ such that 
\begin{eqnarray}\label{eq4.28-1}\|\nabla^{m-\frac n2+\frac 12} u\|_{L^\infty(\Omega)}&\leq & C_\Omega    \sum_{|\alpha|\leq  m-\frac n2+\frac 12}
\| d(\cdot)^{m-\frac n2-\frac 12-|\alpha|}f_{\alpha}\|_{L^p(\Omega)},\,\, 
\end{eqnarray}

\noindent for $p>\frac{n}{n-1},$ provided that the norms on the right-hand side of \eqref{eq4.28-1} are finite. 
\end{proposition}

\bp
 Indeed, the integral representation formula
\begin{equation}\label{eq4.29}
u(x)=\int_{\Omega} G(x,y) \sum_{|\alpha|\leq m-\frac n2+\frac 12}c_\alpha \partial^{\alpha}f_{\alpha}(y)
\,dy,\qquad x\in\Omega,
\end{equation}

\noindent follows directly from the definition of Green function.
It implies that
\begin{eqnarray}\label{eq4.30}
\nabla^{m-\frac n2+\frac 12} u(x)&=& \sum_{|\alpha|\leq m-\frac n2+\frac 12}c_\alpha (-1)^{|\alpha|} \int_{\Omega}\nabla_x^{m-\frac n2+\frac 12}\partial^\alpha_y  G(x,y)f_\alpha(y)\,dy.\end{eqnarray}

Furthermore, due to the estimate \eqref{Gr1} with $i=j=m-\frac n2+\frac 12$ we have
\begin{equation}\label{eq4.29-1}
\int_{\Omega}\Bigl|\nabla_x^{m-\frac n2+\frac 12}\nabla_y^{m-\frac n2+\frac 12}G(x,y)\Bigr|\, |f(y)|\,dy\leq C \int_\Omega \frac{|f(y)|}{|x-y|}\,dy, \end{equation}

\noindent while the bounds in \eqref{Gr2} can be used to show that for every $j\leq m-\frac n2-\frac 12$
\begin{eqnarray}\label{eq4.29-2}\nonumber
&&\int_{\Omega}\Bigl|\nabla_x^{m-\frac n2+\frac 12}\nabla_y^j G(x,y)\Bigr|\, |f(y)|\,dy \leq  C\int_{\Omega} \min\left\{1,\Bigg(\frac{d(y)}{|x-y|}\Bigg)^{m-\frac n2+\frac 12-j}\right\}\times\\[4pt] &&\quad \,\times\frac{1}{|x-y|^{\frac n2-m+\frac 12+j}}\min\left\{\frac{|x-y|}{d(x)},\frac{|x-y|}{d(y)},1\right\}^{\frac n2-m+\frac 12+j}|f(y)|\,dy.
\end{eqnarray}

\noindent We split the latter integral to the cases $|x-y|\leq N^{-1}d(y)$ and $|x-y|\geq N^{-1}d(y)$ with $N\geq 25$ (as in Theorem~\ref{p4.1}). Recall that according to \eqref{eq4.22} in the first case $d(x)\approx d(y)$ and therefore 
\begin{equation}\label{eq4.29-4}
\min\left\{\frac{|x-y|}{d(x)},\frac{|x-y|}{d(y)},1\right\}\approx \frac{|x-y|}{d(y)} \quad\mbox{when}\quad |x-y|\leq N^{-1}d(y),\end{equation}

\noindent while in the second case $d(x)\leq |x-y|+d(y)\leq (1+N)|x-y|$, so that 
\begin{equation}\label{eq4.29-5}
\min\left\{\frac{|x-y|}{d(x)},\frac{|x-y|}{d(y)},1\right\}\approx C \quad\mbox{when}\quad |x-y|\geq N^{-1}d(y).\end{equation}

\noindent Hence, the expression on the right-hand side of \eqref{eq4.29-2} can be further estimated by  
\begin{eqnarray}\label{eq4.29-3}\nonumber
&&C\int_{y\in\Omega:\,|x-y|\leq N^{-1}d(y)} \frac{1}{|x-y|^{\frac n2-m+\frac 12+j}}\left(\frac{|x-y|}{d(y)}\right)^{\frac n2-m+\frac 12+j}|f(y)|\,dy\\[4pt] \nonumber &&\qquad + \,C\int_{y\in\Omega:\,|x-y|\geq N^{-1}d(y)} \Bigg(\frac{d(y)}{|x-y|}\Bigg)^{m-\frac n2+\frac 12-j}\frac{1}{|x-y|^{\frac n2-m+\frac 12+j}}\,|f(y)|\,dy\\[4pt]\nonumber &&\leq C\int_{y\in\Omega:\,|x-y|\leq N^{-1}d(y)} d(y)^{m-\frac n2-\frac 12-j}|f(y)|\,dy \\[4pt]\nonumber &&+ \,C\int_{y\in\Omega:\,|x-y|\geq N^{-1}d(y)} d(y)^{m-\frac n2+\frac 12-j}\,\frac{|f(y)|}{|x-y|}\,dy\\[4pt] && \leq C  \int_\Omega d(y)^{m-\frac n2+\frac 12-j}\,\frac{|f(y)|}{|x-y|}\,dy, 
\end{eqnarray}

\noindent as desired.

This finishes the proof of \eqref{eq4.28} and \eqref{eq4.28-1} follows from it via the mapping properties of the Riesz potential. \ep

Proposition~\ref{pGrAp1} has a natural analogue for the case when the dimension is even. The details are as follows.

\begin{proposition}\label{pGrAp2} Let $\Omega\subset {\mathbb{R}}^n$ be a bounded
domain and assume that  $n\in [2,2m] \cap\NN$ is even.  
Consider the
boundary value problem
\begin{equation}\label{eq4.27-e}
(-\Delta)^m u=\sum_{|\alpha|\leq m-\frac n2}c_\alpha \partial^{\alpha}f_{\alpha}\in W^{-m,2}(\Omega), \quad u\in \ring W^{m,2}(\Omega).
\end{equation}

\noindent Then the solution satisfies the estimate
\begin{equation}\label{eq4.28-e}|\nabla^{m-\frac n2} u(x)|\leq  C  \sum_{|\alpha|\leq  m-\frac n2}
\int_\Omega d(y)^{m-\frac n2-|\alpha|}\,\log\left(1+\frac{d(y)}{|x-y|}\right)|f_\alpha(y)|\,dy,\,  
\end{equation}

\noindent for all $x\in\Omega,$ whenever the integrals on the right-hand side of \eqref{eq4.28-e} are finite. The constant $C$ in \eqref{eq4.28-e} depends on $m$ and $n$ only.

In particular, for every $\eps\in(0,1)$ there exists a constant $C_{\Omega,\eps}>0$ depending on $m$, $n$, $\eps$ and the domain $\Omega$ such that 
\begin{eqnarray}\label{eq4.28-1-e}\|\nabla^{m-\frac n2} u\|_{L^\infty(\Omega)}&\leq & C_{\Omega,\eps} \sum_{|\alpha|\leq  m-\frac n2}
\left\| d(y)^{m-\frac n2-|\alpha|+\eps} f_\alpha\right\|_{L^p(\Omega)},
\end{eqnarray}

\noindent for all $p>\frac{n}{n-\eps}$, provided that the norms on the right-hand side of \eqref{eq4.28-1-e} are finite. 
\end{proposition}

\bp The argument is fairly close to the proof of Proposition~\ref{pGrAp1}. 
We write 
\begin{eqnarray}\label{eq4.30-e}
\left|\nabla^{m-\frac n2} u(x)\right|&\leq & C \sum_{|\alpha|\leq m-\frac n2}\int_{\Omega}\Bigl|\nabla_x^{m-\frac n2}\nabla^{|\alpha|}_y  G(x,y)\Bigr|\, |f_\alpha(y)|\,dy,\end{eqnarray}

\noindent for every $x\in\Omega$, and split the integrals on the right-hand side according to whether $|x-y|\leq N^{-1}d(y)$ or $|x-y|\geq N^{-1}d(y)$, $N\geq 25$. Then using \eqref{eq4.29-4} and \eqref{eq4.29-5} we bound each term on the right-hand side of \eqref{eq4.30-e} by 

\begin{eqnarray}\label{eq4.29-3-e}\nonumber
&&C\int_{y\in\Omega:\,|x-y|\leq N^{-1}d(y)} \frac{1}{|x-y|^{\frac n2-m+|\alpha|}}\left(\frac{|x-y|}{d(y)}\right)^{\frac n2-m+|\alpha|}\times\\[4pt]\nonumber&&\qquad \times
\log\left(1+\frac{d(y)}{|x-y|}\right) |f_\alpha(y)|\,dy\\[4pt] 
&& + \,C\int_{y\in\Omega:\,|x-y|\geq N^{-1}d(y)} \Bigg(\frac{d(y)}{|x-y|}\Bigg)^{m-\frac n2-|\alpha|}\frac{1}{|x-y|^{\frac n2-m+|\alpha|}}\times\nonumber\\[4pt]
&&\qquad\times\log\left(1+\frac{\min\{d(y),d(x)\}}{|x-y|}\right)\, |f_\alpha(y)|\,dy,\end{eqnarray}

\noindent However, if $|x-y|\geq N^{-1}d(y)$ and hence, $d(x)\leq (N+1)|x-y|$, we have 
\begin{equation}\label{eq4.27-e0}
\log\left(1+\frac{\min\{d(y),d(x)\}}{|x-y|}\right)\approx C\approx \log\left(1+\frac{d(y)}{|x-y|}\right).
\end{equation}

\noindent Therefore, both terms in \eqref{eq4.29-3-e}are bounded by 
\begin{eqnarray}\label{eq4.29-3-e0} C\int_{\Omega} d(y)^{m-\frac n2-j}\log\left(1+\frac{d(y)}{|x-y|}\right) |f_\alpha(y)|\,dy , 
\end{eqnarray}

\noindent which leads to \eqref{eq4.28-e}. 

Finally, for every $0<\eps<1$ there is $C_\eps>0$ such that $\log (1+x)\leq C_\eps x^\eps,$ $x>0$, which implies that
\begin{eqnarray}\label{eq4.28-e1}|\nabla^{m-\frac n2} u(x)|&\leq & C_\eps  \sum_{|\alpha|\leq  m-\frac n2}
\int_\Omega d(y)^{m-\frac n2-|\alpha|}\,\left(\frac{d(y)}{|x-y|}\right)^\eps |f_\alpha(y)|\,dy,
\end{eqnarray}

\noindent for all $x\in\Omega,$ $0<\eps<1$.

Then, by the mapping properties of the Riesz potential we recover an estimate 
\begin{equation}\label{eq4.28-e2}\|\nabla^{m-\frac n2} u\|_{L^\infty(\Omega)}\leq  C_{\Omega,\eps}  \sum_{|\alpha|\leq  m-\frac n2}
\left\| d(y)^{m-\frac n2-|\alpha|+\eps} f_\alpha\right\|_{L^p(\Omega)},\,   p>\frac{n}{n-\eps},
\end{equation}

\noindent which leads to \eqref{eq4.28-1-e}.\ep


\end{document}